\documentclass[11pt]{amsart}

\newif\ifdraft\draftfalse

\usepackage{amssymb}
\usepackage{amsthm}
\usepackage{eucal}
\usepackage[english]{babel}
\usepackage{nicefrac}
\usepackage{times}
\usepackage{latexsym,amscd,amsmath,amsthm,amssymb}

\makeatletter
\def\@begintheorem#1#2[#3]{%
    \def\naam{#1}
  \deferred@thm@head{\the\thm@headfont \thm@indent
    \@ifempty{#1}{\let\thmname\@gobble}{\let\thmname\@iden}%
    \@ifempty{#2}{\let\thmnumber\@gobble}{\let\thmnumber\@iden}%
    \@ifempty{#3}{\let\thmnote\@gobble}{\let\thmnote\@iden}%
    \thm@swap\swappedhead\thmhead{#1}{#2}{#3}%
    \the\thm@headpunct
    \thmheadnl 
    \hskip\thm@headsep
  }%
  \ignorespaces}
\makeatother

\newcommand{\kantlijndraft}[1]{\ifdraft\hspace{-\lastskip}%
\vadjust{\vspace{-1mm}\smash{\llap{{\tt #1}\hspace{8mm}}}\vspace{1mm}}\fi}

\def\voegToe#1#2#3{\immediate\write1{\string\newlabel{#1}{{#2}{#3}}}}

\newcommand{\thlabel}[1]{\voegToe{#1}{\naam\noexpand~\thetheorem}{\thepage}\kantlijndraft{#1}}

\makeatletter
\renewcommand{\label}[1]{\voegToe{#1}{\@currentlabel}{\thepage}\kantlijndraft{#1}}
\makeatother

\newtheorem{theorem}{Theorem}[section]
\newtheorem{lemma}[theorem]{Lemma}
\newtheorem{corollary}[theorem]{Corollary}
\newtheorem{question}[theorem]{Question}
\newtheorem{proposition}[theorem]{Proposition}
\theoremstyle{definition}
\newtheorem{observation}[theorem]{Observation}
\newtheorem{example}[theorem]{Example}
\newtheorem{definition}[theorem]{Definition}
\theoremstyle{remark}

\numberwithin{equation}{section}

\newtheorem{claim2}{\sc Claim}



\newcommand{\sse}{\subseteq}						
\newcommand{\minus}{\backslash}						
\newcommand{\Un}{\bigcup}							
\newcommand{\un}{\cup}								
\newcommand{\Meet}{\bigcap}							
\newcommand{\meet}{\cap}							
\newcommand{\es}{\varnothing}						

\newcommand{\cl}[1]{\ensuremath{\overline{#1}}}

\newcommand{\scr}[1]{\ensuremath{\mathcal{#1}}}

\def\juhasz{Juh{\'a}sz}

\begin{document}

\title{On diagonal degrees and star networks}

 \author{Nathan Carlson}\address{Department of Mathematics, California Lutheran University, 60 W. Olsen Rd, MC 3750, 
Thousand Oaks, CA 91360 USA}
\email{ncarlson@callutheran.edu}

\begin{abstract}
Given an open cover $\scr{U}$ of a topological space $X$, we introduce the notion of a star network for $\scr{U}$. The associated cardinal function $sn(X)$, where $e(X)\leq sn(X)\leq L(X)$, is used to establish new cardinal inequalities involving diagonal degrees. We show $|X|\leq sn(X)^{\Delta(X)}$ for a $T_1$ space $X$, giving a partial answer to a long-standing question of Angelo Bella. Many further results are given using variations of $sn(X)$. One result has as corollaries Buzyakova's theorem that a ccc space with a regular $G_\delta$-diagonal has cardinality at most $\mathfrak{c}$, as well as three results of Gotchev. Further results lead to logical improvements of theorems of Basile, Bella, and Ridderbos, a partial solution to a question of the same authors, and a theorem of Gotchev, Tkachenko, and Tkachuk. Finally, we define the Urysohn extent $Ue(X)$ with the property $Ue(X)\leq\min\{aL(X),e(X)\}$ and use the Erd\H{o}s-Rado theorem to show that $|X|\leq 2^{Ue(X)\overline{\Delta}(X)}$ for any Urysohn space $X$.
\end{abstract}



\maketitle

\section{Introduction.}

The \emph{diagonal} of a topological space $X$, denoted by $\Delta_X$, is defined as $\Delta_X=\{(x,x)\in X^2:x\in X\}$. Considerable progress has been made in the theory of cardinal inequalities involving the \emph{diagonal degree} $\Delta(X)$ (Definition~\ref{diagdegree}) and its variants. (See~\cite{BBR11},~\cite{Bel87},~\cite{BCG23},~\cite{BS20},~\cite{Buz06},~\cite{Got19}, and others). An early result was given in 1977 by Ginsburg and Woods~\cite{GW77} who used the Erd\H{o}s-Rado theorem (Theorem~\ref{ER}) to demonstrate that if $X$ is $T_1$ then $|X|\leq 2^{e(X)\Delta(X)}$. Here $e(X)$, the \emph{extent} of $X$, is the supremum of the cardinalities of closed discrete subspaces of $X$. In 1989 Bella~\cite{Bel89} asked if $e(X)$ could be moved out of the exponent in this result.


\begin{question}[Bella~\cite{Bel89}]\label{BellaQuestion}
If $X$ is $T_1$ is $|X|\leq e(X)^{\Delta(X)}$?
\end{question}

In 1982 Alas~\cite{Ala82} showed that if $X$ is $T_1$ then $|X|\leq L(X)^{\Delta(X)}$. As  $e(X)\leq L(X)$ for any space $X$, this represents a partial answer to the question of Bella. Proofs of this result also occur in~\cite{Bel89} and~\cite{Xus19}. In this study we use the notion of a \emph{star network} (Definition~\ref{starnetwork}) to define the cardinal function $sn(X)$ with the property $e(X)\leq sn(X)\leq L(X)$ and show if $X$ is $T_1$ then $|X|\leq sn(X)^{\Delta(X)}$ (Theorem~\ref{snbound}). This gives an improved partial solution to the question of Bella. 

For each integer $m\geq 1$ we define the cardinal function $sn_m(X)$ using \emph{star}-$m$ \emph{networks} (Definition~\ref{starnetwork}). There is a fundamental  ``pairing" of each $sn_m(X)$ with each corresponding rank $m$ diagonal degree $\Delta_m(X)$ (Definition~\ref{rank}) in the sense that $|X|\leq sn_m(X)^{\Delta_m(X)}$ for any space $X$ for which $\Delta_m(X)$ is defined (Theorem~\ref{general}). Likewise, \emph{weak star}-$m$ \emph{networks} (Definition~\ref{starnetwork}) are used to define $wsn_m(X)$, which is ``paired" with the strong rank $m$ diagonal degree $s\Delta_m(X)$ through the cardinal inequality $|X|\leq wsn_m(X)^{s\Delta_m(X)}$ for any space $X$ for which $s\Delta_m(X)$ is defined (Theorem~\ref{generaltheta}). 

Many cardinal inequalities are established using $sn_m(X)$ and $wsn_m(X)$ in~\S3. For example, we show $sn_2(X)\leq we(X)$ (Theorem~\ref{inequalitystring1}), where $we(X)$ is the \emph{weak extent} of $X$ (Definition~\ref{weakextent}). Therefore the result that if $X$ is Hausdorff then $|X|\leq sn_2(X)^{\Delta_2(X)}$ is a logical improvement of the result of Basile, Bella, and Ridderbos~\cite{BBR11} that the cardinality of a Hausdorff space is bounded by $we(X)^{\Delta_2(X)}$. Another example is Theorem~\ref{sn_2thetabound}: $wsn_2(X)\leq 2^{wL(X)}$ for any space $X$. A consequence is that if $X$ is Urysohn then $|X|\leq 2^{wL(X)s\Delta_2(X)}$. This gives a partial answer to Question 4.8 in~\cite{BBR11}. A third example is Theorem~\ref{somewsn2bound}: $wsn_2(X)\leq wL(X)^{dot(X)}$, which has as a corollary the result of Gotchev, Tkachenko, and Tkachuk~\cite{GTT16} that if $X$ is Urysohn then $|X|\leq wL(X)^{dot(X)s\Delta_2(X)}$.

In Definition~\ref{Uryextent} we define the \emph{Urysohn extent} $Ue(X)$ with the property $Ue(X)\leq\min\{aL(X),e(X)\}$ and use the Erd\H{o}s-Rado theorem to show that if $X$ is Urysohn then $|X|\leq 2^{Ue(X)\overline{\Delta}(X)}$ (Theorem~\ref{Ue(X)}). This appears to be the first application of the Erd\H{o}s-Rado theorem in connection with the regular diagonal degree $\overline{\Delta}(X)$.

In~\S4 we use \emph{regular star networks} (Definition~\ref{regularstarnetwork}) to define the cardinal function $\overline{sn}(X)$ and show there is also a natural ``pairing" of $\overline{sn}(X)$ with $\overline{\Delta}(X)$. This is demonstrated in Theorem~\ref{overlinesnCardBound}: If $X$ is Urysohn then $|X|\leq\overline{sn}(X)^{\overline{\Delta}(X)}$. By utilizing this result and the fact that $\overline{sn}(X)\leq 2^{c(X)}$ for any space $X$, we have as a corollary Buzyakova's theorem that a ccc space with a regular $G_\delta$-diagonal has cardinality at most $\mathfrak{c}$. We also have as corollaries three results of Gotchev~\cite{Got19}: If $X$ is Urysohn then $|X|\leq aL(X)^{\overline{\Delta}(X)}$, $|X|\leq wL(X)^{\chi(X)\overline{\Delta}(X)}$, and $|X|\leq 2^{\overline{\Delta}(X)2^{wL(X)}}$.

In this paper we make no global assumption on any separation axiom on a topological space.

\section{Definitions.}

In this section we give the main definitions used in this paper, as well as several characterizations. For all other notions not defined here, see Engelking~\cite{Engelking} and~\juhasz~\cite{Juh80}.

\begin{definition}\label{diagdegree}
The \emph{diagonal} of a topological space $X$ is defined to be $\Delta_X=\{(x,x)\in X^2:x\in X\}$. For a cardinal $\kappa$, $X$ has a $G_\kappa$-\emph{diagonal} if there exists a family $\{U_\alpha:\alpha<\kappa\}$ of open sets in $X^2$ such that $\Delta_X=\Meet_{\alpha<\kappa}U_\alpha$. If additionally $\Delta_X=\Meet_{\alpha<\kappa}U_\alpha=\Meet_{\alpha<\kappa}\cl{U_\alpha}$ then $X$ has a \emph{regular} $G_\kappa$-\emph{diagonal}. The \emph{diagonal degree of }$X$, denoted by $\Delta(X)$, is the least infinite cardinal $\kappa$ such that $X$ has a $G_\kappa$-diagonal. The \emph{regular diagonal degree of }$X$, denoted by $\overline{\Delta}(X)$, is the least infinite cardinal $\kappa$ such that $X$ has a regular $G_\kappa$-diagonal. 
\end{definition}

It is not hard to see that a space $X$ has a $G_\kappa$-diagonal for some infinite cardinal $\kappa$ if and only if $X$ is $T_1$. Likewise, $X$ has a regular $G_\kappa$-diagonal for some infinite cardinal $\kappa$ if and only if $X$ is Urysohn. It follows that $\Delta(X)$ is defined if and only if $X$ is $T_1$ and that $\overline{\Delta}(X)$ is defined if and only if $X$ is Urysohn.

\begin{definition}\label{star}
Let $X$ be a space, $\scr{U}$ a family of subsets of $X$, and $A\sse X$. We define $St(A,\scr{U})=\Un\{U\in\scr{U}:U\meet A\neq\es\}$ and if $A=\{a\}$ we rewrite $St(\{a\},\scr{U})=St(a,\scr{U})$. We define $St^2(A,\scr{U})=St(St(A,\scr{U}),\scr{U})$ and, for an integer $n>2$, we define $St^n(A,\scr{U})=St(St^{n-1}(A,\scr{U}),\scr{U})$. Last, we define $St^0(A,\scr{U})=A$ and $St^1(A,\scr{U})=St(A,\scr{U})$.
\end{definition}

Bonanzinga~\cite{Bon98} defined the following cardinal function.

\begin{definition}\label{stL}
Given a space $X$, $stL(X)$ is defined to be the least infinite cardinal $\kappa$ such that for every open cover $\scr{U}$ there exists $A\sse [X]^{\leq\kappa}$ such that $St(A,\scr{U})= X$.
\end{definition}

It is not hard to see that $stL(X)$ is also the least infinite cardinal $\kappa$ such that for every open cover $\scr{U}$ there exists $A\sse [X]^{\leq\kappa}$ such that $\{St(x,\scr{U}):x\in A\}$ covers $X$. Observe that $stL(X)\leq L(X)$.

\begin{definition}\label{rank}
Let $X$ be a space, let $n\in\omega$, and let $\kappa$ be an infinite cardinal. We say $X$ has a \emph{rank }$n$ $G_\kappa$-\emph{diagonal} (\emph{strong rank }$n$ $G_\kappa$-\emph{diagonal}) if there is a family $\{\scr{U}_\alpha:\alpha<\kappa\}$ of open covers of $X$ such that for all $x\neq y$ there exists $\alpha<\kappa$ such that $y\notin St^n(x,\scr{U}_\alpha)$ ($y\notin\cl{St^n(x,\scr{U}_\alpha)}$). If $\kappa=\omega$ then we denote a rank $n$ $G_\delta$-diagonal by \emph{rank }$n$ \emph{diagonal}. The least infinite cardinal $\kappa$ such that $X$ has a rank $n$ $G_\kappa$-diagonal or a strong rank $n$ $G_\kappa$-diagonal is denoted by $\Delta_n(X)$ and $s\Delta_n(X)$, respectively. 
\end{definition}

It is clear that $\Delta_n(X)\leq\min\{\Delta_{n+1}(X),s\Delta_n(X)\}$. We have the following characterization. It is a straightforward generalization of a result of Ceder.

\begin{proposition}[Ceder~\cite{Ced61}]\label{ceder}
A space $X$ has a $G_\kappa$-diagonal if and only if it has a rank 1 $G_\kappa$-diagonal.
\end{proposition}

Another way to state Proposition~\ref{ceder} is that a space $X$ has a $G_\kappa$-diagonal if and only if there is a family of open covers $\{\scr{U}_\alpha:\alpha<\kappa\}$ such that $\{x\}=\Meet_{\alpha<\kappa}St(x,\scr{U}_\alpha)$ for every $x\in X$. It follows that $\psi(X)\leq\Delta(X)=\Delta_1(X)$ for any $T_1$ space $X$. Additionally, one can see that $\psi_c(X)\leq s\Delta(X)$ if $X$ is Hausdorff, where $s\Delta(X)$ is defined to be $s\Delta_1(X)$. This is because if $\{\scr{U}_\alpha:\alpha<\kappa\}$ is a family of open covers witnessing that $s\Delta(X)=\kappa$, then for all $x\in X$ we have that $\{St(x,\scr{U}_\alpha):\alpha<\kappa\}$ forms a closed pseudobase at $x$.

For a space $X$ and a cardinal $\kappa$, we have the following two similar characterizations of the statement ``$X$ has a regular $G_\kappa$-diagonal" in Proposition~\ref{gotchev}. (Notice that $U$ and $V$ are in $\scr{U}_\alpha$ in the second characterization and not necessarily in $\scr{U}_\alpha$ in the first). Zenor proved the first in~\cite{Zen72} in the case $\kappa=\omega$, which easily generalizes, and Gotchev proved the second in~\cite{Got19}. It is interesting that the second characterization is logically stronger than the first but both are still equivalent. 

\begin{proposition}[Zenor~\cite{Zen72} and Gotchev~\cite{Got19}]\label{gotchev}
Let $X$ be a space and let $\kappa$ be an infinite cardinal. The following are equivalent.
\begin{enumerate}
\item $X$ has a regular $G_\kappa$-diagonal.
\item There is a family $\{\scr{U}_\alpha:\alpha<\kappa\}$ of open covers of $X$ such that if $x\neq y\in X$ then there exists $\alpha<\kappa$ and open sets $U$ and $V$ containing $x$ and $y$ respectively such that $V\meet\cl{St(U,\scr{U}_\alpha)}=\es$.
\item There is a family $\{\scr{U}_\alpha:\alpha<\kappa\}$ of open covers of $X$ such that if $x\neq y\in X$ then there exists $\alpha<\kappa$ and open sets $U$ and $V$ in $\scr{U}_\alpha$ containing $x$ and $y$ respectively such that $V\meet\cl{St(U,\scr{U}_\alpha)}=\es$.
\end{enumerate}
\end{proposition}

It follows by Proposition~\ref{gotchev} that $s\Delta(X)\leq\overline{\Delta}(X)\leq s\Delta_2(X)$ if $X$ is Urysohn. Therefore if a space has strong rank 2-diagonal then it has a regular $G_\delta$-diagonal. However, the following difficult question of Bella remains open.

\begin{question}[Bella~\cite{Bel87}]
Does any space with a regular $G_\delta$-diagonal have a rank 2-diagonal?
\end{question}

\begin{definition}\label{starnetwork}
Given a cover $\scr{U}$ of a space $X$ and an integer $m\geq 1$, we define a \emph{star}-$m$\emph{ network} for $\scr{U}$ to be a collection $\scr{N}$ of subsets of $X$ such that for all $x\in X$ there exists $N\in\scr{N}$ such that $x\in N\sse St^m(x,\scr{U})$. We define a \emph{weak star-}$m$ \emph{network} for $\scr{U}$ in a similar way except $x\in N\sse\cl{St^m(x,\scr{U})}$. By \emph{star network} we mean a star-$1$ network and by \emph{weak star network} we mean a weak star-$1$ network. We define $sn_m(X)$ $(wsn_m(X))$ to be the least infinite cardinal $\kappa$ such that every open cover $\scr{U}$ has a star-$m$ network (weak star-$m$ network) of cardinality at most $\kappa$. We denote $sn_1(X)$ by $sn(X)$ and $wsn_1(X)$ by $wsn(X)$.
\end{definition}

The term ``star network" is used above because of its similarity to a network for a space. A network $\scr{N}$ has the property that whenever $x\in U$ for an open set $U$ there exists $N\in\scr{N}$ such that $x\in N\sse U$. Similarly, a star network $\scr{N}$ for an open cover $\scr{U}$ has the property that for all $x\in X$ there exists an $N\in\scr{N}$ such that $x\in N\sse St(x,\scr{U})$. While the terminology ``star $\scr{P}$" has a particular meaning when $\scr{P}$ is a property of a space $X$, we remark that ``network" is not a property of a space and so should not be confused with this terminology. 

For an integer $m\geq 1$, every star-$m$ network is a star-$(m+1)$ network and therefore $sn_{m+1}(X)\leq sn_m(X)$. In particular $sn_2(X)\leq sn(X)$. Likewise, $wsn_{m+1}(X)\leq wsn_m(X)$ and thus $wsn_2(X)\leq wsn(X)$. Observe also that $wsn_m(X)\leq sn_m(X)$ for any space $X$, in particular $wsn(X)\leq sn(X)$.

\begin{definition}\label{regularstarnetwork}
Given a cover $\scr{U}$ of a space $X$, a \emph{regular star network for} $\scr{U}$ is a collection $\scr{N}$ of subsets of $X$ such that for all $x\in X$ there exists $N\in\scr{N}$ such that $x\in N\sse\Meet\{\cl{St(U,\scr{U})}:x\in U\in\scr{U}\}$. Then $\overline{sn}(X)$ is defined to be the least infinite cardinal $\kappa$ such that every open cover $\scr{U}$ has a regular star network of cardinality at most $\kappa$.
\end{definition}

Note that $wsn_2(X)\leq\overline{sn}(X)\leq wsn(X)$ for any space $X$. This is because $\cl{St(x,\scr{U})}\sse\Meet\{\cl{St(U,\scr{U})}:x\in U\in\scr{U}\}\sse\cl{St^2(x,\scr{U})}$.

Recall that the \emph{extent} $e(X)$ of a space $X$ is defined as the supremum of the cardinalities of the closed discrete subsets of $X$. We define the \emph{Urysohn extent} in Definition~\ref{Uryextent}. It uses the notion of a Urysohn discrete set defined by Schr\"oder in~\cite{Sch93}. First recall the following.

\begin{definition}\label{thetaacc}
Given a set $A$ in a space $X$, the $\theta$-\emph{closure} of $A$ is defined as $cl_\theta(A)=\{x\in X:\cl{U}\meet A\neq\es\textup{ whenever }U\textup{ is an open neighborhood of }x\}$. $A$ is $\theta$-\emph{closed} if $A=cl_\theta(A)$. A point $x$ is a $\theta$-\emph{accumulation point} of $A$ if $\cl{U}\meet A\minus\{x\}\neq\es$ for every open set $U$ containing $x$.
\end{definition}

\begin{definition}[Schr\"oder~\cite{Sch93}]\label{Urydiscrete}
A set $D$ in a space $X$ is \emph{Urysohn discrete} if for all $d\in D$ there exists an open set $U_d$ containing $d$ such that $\cl{U_d}\meet (D\minus\{d\})=\es$.
\end{definition}

\begin{definition}\label{Uryextent}
Given a space $X$, we define the \emph{Urysohn extent of }$X$ by 
$$Ue(X)=\sup\{|D|: D\textup{ is }\theta\textup{-closed and Urysohn discrete}\}.$$ 
\end{definition}

One can see that $Ue(X)\leq e(X)$. Our last variant on the extent is the \emph{weak extent} $we(X)$ of a space $X$.

\begin{definition}\label{weakextent}
The \emph{weak extent} of a space $X$, denoted by $we(X)$, is the least cardinal $\kappa$ such that for every open cover $\scr{U}$ of $X$ there is a subset $A$ of $X$ of cardinality no greater than $\kappa$ such that $St(A, \scr{U}) = X$. 
\end{definition}

It is clear that $we(X)\leq d(X)$ and $we(X)\leq e(X)$. Note that spaces with countable weak extent are called \emph{star countable} by some authors. 

Recall the following definitions of $wL(X)$, $aL(X)$, and $wL_c(X)$.

\begin{definition}\label{wLaL}
Given a space $X$, the \emph{weak Lindel\"of} degree $wL(X)$ is defined to be the least infinite cardinal $\kappa$ such that every open cover $\scr{U}$ of $X$ has a subfamily $\scr{V}$ of cardinality at most $\kappa$ such that $X=\cl{\Un\scr{V}}$. The \emph{almost Lindel\"of} degree $aL(X)$ is defined to be the least infinite cardinal $\kappa$ such that every open cover $\scr{U}$ of $X$ has a subfamily $\scr{V}$ of cardinality at most $\kappa$ such that $X=\Un_{V\in\scr{V}}\cl{V}$.
\end{definition}

\begin{definition}\label{wLcdef}
Given a space $X$ and a subset $A$ of $X$, we define $wL(A,X)$ to be the least infinite cardinal $\kappa$ such that every open cover $\scr{U}$ of $A$ has a subfamily $\scr{V}$ of cardinality at most $\kappa$ such that $A\sse\cl{\Un\scr{V}}$. The \emph{weak Lindel\"of degree of} $X$ \emph{with respect to closed sets} is $wL_c(X)=\sup\{wL(C,X):C\textup{ is a closed subset of }X\}$.
\end{definition}

It is clear that $wL(X)\leq aL(X)\leq L(X)$ and $wL(X)\leq wL_c(X)\leq c(X)$. 

\begin{definition}[\cite{C23}]\label{wpsi}
Let $X$ be a Hausdorff space and let $x\in X$. A collection of open sets $\scr{V}$ is a \emph{weak closed pseudobase at }$x$ if $\{x\}=\Meet_{V\in\scr{V}}\cl{V}$. (Note that it is not necessarily the case that $x\in V$ for any $V\in\scr{V}$). We define $w\psi_c(x,X)$ to be the least infinite cardinal $\kappa$ such that $x$ has a weak closed pseudobase of cardinality $\kappa$. The \emph{weak closed pseudocharacter} $w\psi_c(X)$ is defined as $w\psi_c(X)=\sup\{w\psi_c(x,X):x\in X\}$.
\end{definition}

It is clear that $w\psi_c(X)\leq\psi_c(X)$ for any Hausdorff space $X$. The following was defined in~\cite{GTT16}. It can be seen that $dot(X)\leq\pi\chi(X)$ and $dot(X)\leq c(X)$.

\begin{definition}[\cite{GTT16}]\label{dot}
Let $X$ be a space. The \emph{dense o-tightness} of $X$, denoted by $dot(X)$, is the least infinite cardinal $\kappa$ such that whenever $x\in X=\cl{\Un\scr{U}}$, for an open family $\scr{U}$, there exists $\scr{V}\sse\scr{U}$ such that $|\scr{V}|\leq\kappa$ and $x\in\cl{\Un\scr{V}}$.
\end{definition}

\section{On star networks, weak star networks, and $Ue(X)$.}

In this section we establish several cardinal inequalities involving $sn(X)$ and its variants, including a partial answer to Question~\ref{BellaQuestion} of Bella.  We begin with preliminary inequalities.

\begin{proposition}\label{basicsn}
$e(X)\leq sn(X)\leq L(X)$ for any space $X$.
\end{proposition}

\begin{proof}
First we show $sn(X)\leq L(X)$. Let $\kappa=L(X)$ and let $\scr{U}$ be an open cover of $X$. As $\kappa=L(X)$ there exists $\scr{V}\sse\scr{U}$ such that $|\scr{V}|\leq\kappa$ and $X=\Un\scr{V}$. If $x\in X$ then there exists $V\sse\scr{V}$ such that $x\in V\sse St(x,\scr{U})$. This shows $\scr{V}$ is a star network for $\scr{U}$ and that $sn(X)\leq |\scr{V}|\leq\kappa$.

To show $e(X)\leq sn(X)$, let $\kappa=sn(X)$ and suppose by way of contradiction that $e(X)\geq\kappa^+$. Then there exists a closed discrete set $D$ such that $|D|=\kappa^+$. For all $d\in D$ there exists an open set $U_d$ containing $d$ such that $U_d\meet D\minus\{d\}=\es$. Then $\scr{U}=\{X\minus D\}\un\{U_d:d\in D\}$ is an open cover of $X$. Observe that $St(d,\scr{U})=U_d$ for all $d\in D$.

As $sn(X)=\kappa$ there exists a star network $\scr{N}$ for $\scr{U}$ such that $|\scr{N}|\leq\kappa$. It follows that for all $d\in D$ there exists $N_d\in\scr{N}$ such that $d\in N_d\sse St(d,\scr{U})=U_d$. As $|D|=\kappa^+$ and $|\scr{N}|\leq\kappa$, there exists $d_1\neq d_2\in D$ such that $N_{d_1}=N_{d_2}$. Let $N=N_{d_1}=N_{d_2}$. Then $\{d_1,d_2\}\sse N\sse U_{d_1}\meet U_{d_2}$. This implies that $d_2\in U_{d_1}$, a contradiction. Therefore $e(X)\leq\kappa=sn(X)$.
\end{proof}

The next theorem gives a partial solution to Question~\ref{BellaQuestion} of Bella. The proof is modeled after the proof that $|X|\leq nw(X)^{\psi(X)}$ for a $T_1$ space $X$, where $\{St(x,\scr{U}_\alpha):\alpha<\kappa\}$ is a used as a pseudobase at $x\in X$. (See~\cite{Juh80}).

\begin{theorem}\label{snbound}
If $X$ is $T_1$ then $|X|\leq sn(X)^{\Delta(X)}$.
\end{theorem}

\begin{proof}
Let $\kappa=\Delta(X)$ and let $\lambda=sn(X)$. By Proposition~\ref{ceder}, there exist open covers $\{\scr{U}_\alpha:\alpha<\kappa\}$ of $X$ such that for all $x\neq y\in X$ there exists an $\alpha<\kappa$ such that $y\notin St(x,\scr{U}_\alpha)$. As $\lambda=sn(X)$, for all $\alpha<\kappa$ there exists a star network $\scr{N}_\alpha$ for $\scr{U}_\alpha$ such that $|\scr{N}_\alpha|\leq\lambda$. Let $\scr{N}=\Un_{\alpha<\kappa}\scr{N}_\alpha$ and note $|\scr{N}|\leq\lambda\kappa$. Let $\scr{C}=\{\Meet\scr{M}:\scr{M}\in[\scr{N}]^{\leq\kappa}\}$. Then $|\scr{C}|\leq |\scr{N}|^\kappa\leq(\lambda\kappa)^\kappa=\lambda^\kappa$.

Fix $x\in X$. For all $\alpha<\kappa$ there exists $N_\alpha\in\scr{N}_\alpha\sse\scr{N}$ such that $x\in N_\alpha\sse St(x,\scr{U}_\alpha)$. It follows that $x\in\Meet_{\alpha<\kappa}N_\alpha\sse\Meet_{\alpha<\kappa}St(x,\scr{U}_\alpha)$. If $y\neq x$ then there exists $\alpha<\kappa$ such that $y\notin St(x,\scr{U}_\alpha)$ and so $y\notin\Meet_{\alpha<\kappa}St(x,\scr{U}_\alpha)$. Therefore $y\notin\Meet_{\alpha<\kappa}N_\alpha$ and $\{x\}=\Meet_{\alpha<\kappa}N_\alpha$. As $\Meet_{\alpha<\kappa}N_\alpha\in\scr{C}$, we have $|X|\leq|\scr{C}|\leq\lambda^\kappa$.
\end{proof}

In light of Proposition~\ref{basicsn}, compare the above result with the following theorem of Ginsburg and Woods.

\begin{theorem}[Ginsburg and Woods~\cite{GW77}]
 If $X$ is $T_1$ then $|X|\leq 2^{e(X)\Delta(X)}$. 
\end{theorem}

Also, it is straightforward to modify the proof of Theorem~\ref{snbound} to show the following.

\begin{theorem}
If $X$ is Hausdorff then $|X|\leq wsn(X)^{s\Delta(X)}$.
\end{theorem}

By Proposition~\ref{basicsn} and Theorem~\ref{snbound}, we have the following result of Xuan and Song.

\begin{corollary}[Xuan and Song~\cite{Xus19}] If $X$ is $T_1$ then $|X|\leq L(X)^{\Delta(X)}$.
\end{corollary}

\begin{question} Is there an example of a $T_1$ space for which $sn(X)^{\Delta(X)}<L(X)^{\Delta(X)}$?
\end{question}

The next theorem establishes nice relationships between $sn_2(X)$, $stL(X)$, and $we(X)$. (See Definitions~\ref{starnetwork},~\ref{stL}, and~\ref{weakextent}, respectively). 

\begin{theorem}\label{inequalitystring1}
$sn_2(X)\leq stL(X)\leq we(X)$ for any space $X$.
\end{theorem}

\begin{proof}
For the first inequality, let $\kappa=stL(X)$ and let $\scr{U}$ be an open cover of $X$. As $stL(X)=\kappa$ there exists $A\in [X]^{\leq\kappa}$ such that $\{St(x,\scr{U}):x\in A\}$ covers $X$. Let $\scr{N}=\{St(x,\scr{U}):x\in A\}$ and note $|\scr{N}|\leq\kappa$. Fix $x\in X$. Then $x\in St(y,\scr{U})$ for some $y\in A$. It follows that there exists $U\in\scr{U}$ such that $x\in U$ and $y\in U$. We want to show that $St(y,\scr{U})\sse St(St(x,\scr{U}),\scr{U})$. Let $W\in\scr{U}$ such that $y\in W$. Then $y\in W\meet U$ and so $W\meet U\neq\es$. Therefore $W\sse St(St(x,\scr{U}),\scr{U})$. This shows $St(y,\scr{U})\sse St(St(x,\scr{U}),\scr{U})$ and that $sn_2(X)\leq stL(X)$.

For the second inequality, let $we(X)=\kappa$ and let $\scr{U}$ be an open cover of $X$. As $we(X)=\kappa$, there exists $A\in [X]^{\leq\kappa}$ such that $St(A,\scr{U})=X$. We show $\{St(a,\scr{U}):a\in A\}$ covers $X$. Let $x\in X$. Then $x\in St(A,\scr{U})$ and there exists $U\in\scr{U}$ such that $x\in U$ and $U\meet A\neq\es$. Let $a\in U\meet A$. Then $x\in U\sse St(a,\scr{U})$. This shows $\{St(a,\scr{U}):a\in A\}$ covers $X$ and that $stL(X)\leq we(X)$.
\end{proof}

In~\cite{BBR11}, Basile, Bella, and Ridderbos showed that if $X$ is Hausdorff then $|X|\leq we(X)^{\Delta_2(X)}$. We give a logical improvement of this result using the cardinal function $sn_2(X)$. 

\begin{theorem}\label{sn2bound}
If $X$ is Hausdorff then $|X|\leq sn_2(X)^{\Delta_2(X)}$.
\end{theorem}

\begin{proof}
Let $\kappa=\Delta_2(X)$ and let $\lambda=sn_2(X)$. As $\Delta_2(X)=\kappa$, there exist open covers $\{\scr{U}_\alpha:\alpha<\kappa\}$ of $X$ such that for all $x\neq y\in X$ there exists an $\alpha<\kappa$ such that $y\notin St^2(x,\scr{U}_\alpha)$. As $\lambda=sn_2(X)$, for all $\alpha<\kappa$ there exists a star-$2$ network $\scr{N}_\alpha$ for $\scr{U}_\alpha$ such that $|\scr{N}_\alpha|\leq\lambda$. Let $\scr{N}=\Un_{\alpha<\kappa}\scr{N}_\alpha$ and note $|\scr{N}|\leq\lambda\kappa$. Let $\scr{C}=\{\Meet\scr{M}:\scr{M}\in[\scr{N}]^{\leq\kappa}\}$. Then $|\scr{C}|\leq |\scr{N}|^\kappa\leq(\lambda\kappa)^\kappa=\lambda^\kappa$.

Fix $x\in X$. For all $\alpha<\kappa$ there exists $N_\alpha\in\scr{N}_\alpha\sse\scr{N}$ such that $x\in N_\alpha\sse St^2(x,\scr{U}_\alpha)$. It follows that $x\in\Meet_{\alpha<\kappa}N_\alpha\sse\Meet_{\alpha<\kappa}St^2(x,\scr{U}_\alpha)$. If $y\neq x$ then there exists $\alpha<\kappa$ such that $y\notin St^2(x,\scr{U}_\alpha)$ and so $y\notin\Meet_{\alpha<\kappa}St^2(x,\scr{U}_\alpha)$. Therefore $y\notin\Meet_{\alpha<\kappa}N_\alpha$ and $\{x\}=\Meet_{\alpha<\kappa}N_\alpha$. As $\Meet_{\alpha<\kappa}N_\alpha\in\scr{C}$, we have $|X|\leq|\scr{C}|\leq\lambda^\kappa$.
\end{proof}

Generalizing the above proof we have the following.

\begin{theorem}\label{general}
Let $m$ be an integer such that $m\geq 1$. If $X$ is a space for which $\Delta_m(X)$ is defined, then $|X|\leq sn_m(X)^{\Delta_m(X)}$.
\end{theorem}

Theorems~\ref{snbound} and~\ref{sn2bound} are corollaries of Theorem~\ref{general}. In addition, straightforward modification of the proof Theorem~\ref{sn2bound} gives the following. 

\begin{theorem}\label{Xboundsn_2theta}
If $X$ is Urysohn then $|X|\leq wsn_2(X)^{s\Delta_2(X)}$.
\end{theorem}

This in turn generalizes to the following.

\begin{theorem}\label{generaltheta}
Let $m$ be an integer such that $m\geq 1$. If $X$ is a space for which $s\Delta_m(X)$ is defined, then $|X|\leq wsn_m(X)^{s\Delta_m(X)}$.
\end{theorem}

\begin{corollary}[Basile, Bella, and Ridderbos~\cite{BBR11}]
If $X$ is Hausdorff then $|X|\leq we(X)^{\Delta_2(X)}$.
\end{corollary}

\begin{proof} By Theorems~\ref{sn2bound} and~\ref{inequalitystring1}, we have
$|X|\leq sn_2(X)^{\Delta_2(X)}\leq we(X)^{\Delta_2(X)}$.
\end{proof}

\begin{theorem}\label{sn_2thetabound}
$wsn_2(X)\leq 2^{wL(X)}$ for any space $X$.
\end{theorem}

\begin{proof}
Let $\kappa=wL(X)$ and let $\scr{U}$ be an open cover of $X$. Then there exists $\scr{V}\sse\scr{U}$ such that $|\scr{V}|\leq\kappa$ and $X=\cl{\Un\scr{V}}$. For all $U\in\scr{U}$ let $\scr{V}_U=\{V\in\scr{V}:V\meet U\neq\es$. As $X=\cl{\Un\scr{V}}$ we see that $\scr{V}_U\neq\es$ for all $U\in\scr{U}$.

We wish to find a collection $\scr{N}$ of subsets of $X$ such that $|\scr{N}|\leq 2^\kappa$ and for all $x\in X$ there exists $N\in\scr{N}$ such that $x\in N\sse\cl{St^2(x,\scr{U})}$. We show that $\scr{N}=\{\cl{\Un\scr{V}_U}:U\in\scr{U}\}$ has these properties. First, note that $\scr{N}\sse\{\cl{\Un\scr{W}}:\scr{W}\in\scr{P}(\scr{V})\}$ and so $|\scr{N}|\leq |\scr{P}(\scr{V})|\leq 2^\kappa$.

Fix $x\in X$ and fix $U\in\scr{U}$ such that $x\in U$. For an open set $W$ containing $x$, we have $x\in W\meet U$ and therefore there exists $V\in\scr{V}$ such that $W\meet U\meet V\neq\es$. It follows that $V\in\scr{V}_U$. Since $W\meet V\neq\es$, we have $W\meet\Un\scr{V}_U\neq\es$. Therefore $x\in\cl{\Un\scr{V}_U}$.

We show $\Un\scr{V}_U\sse\cl{St^2(x,\scr{U})}$. Let $y\in\Un\scr{V}_U$. Then there exists $V\in\scr{V}_U$ such that $y\in V$ and $U\meet V\neq\es$. Let $T$ be an open set containing $y$. We show $T\meet St^2(x,\scr{U})\neq\es$. As $y\in V\meet T$ we have that $V\meet T\neq\es$. Now as $V\meet U\neq\es$ and $U\sse St(x,\scr{U})$, we have that $V\sse St(St(x,\scr{U}),\scr{U})=St^2(x,\scr{U})$. As $T\meet V\neq\es$, we have that  $T\meet St^2(x,\scr{U})\neq\es$ and $y\in\cl{St^2(x,\scr{U})}$. This shows $\Un\scr{V}_U\sse\cl{St^2(x,\scr{U})}$.

It follows that $x\in\cl{\Un\scr{V}_U}\sse\cl{St^2(x,\scr{U})}$ and that $\scr{N}$ is a weak star-2 network for $\scr{U}$. Therefore $wsn_2(X)\leq |\scr{N}|\leq 2^\kappa=2^{wL(X)}$.
\end{proof}

By Theorems~\ref{sn_2thetabound} and~\ref{Xboundsn_2theta}, we have the following corollary. Corollary~\ref{somecorollary} provides a partial solution to the question of Basile, Bella, and Ridderbos of whether a space $X$ with a strong rank 2 diagonal has cardinality at most $2^{dc(X)}$. (Question 4.8 in~\cite{BBR11}).

\begin{corollary}\label{somecorollary}
If $X$ is Urysohn then $|X|\leq 2^{wL(X)s\Delta_2(X)}$.
\end{corollary}

The proof of the following is a reduced version of the proof of Theorem 4.3 in~\cite{GTT16}. Recall the definition of $dot(X)$ in~\ref{dot}.

\begin{theorem}\label{somewsn2bound}
$wsn_2(X)\leq wL(X)^{dot(X)}$ for any space $X$.
\end{theorem}

\begin{proof}
Let $\kappa=dot(X)$, $\lambda=wL(X)$, and let $\scr{U}$ be an open cover of $X$. As $wL(X)=\lambda$, there exists $\scr{W}\sse\scr{U}$ such that $|\scr{W}|\leq\lambda$ and $X=\cl{\Un\scr{W}}$. As $dot(X)=\kappa$, for all $x\in X$ there exists $\scr{V}_x\in\scr{W}$ such that $x\in\cl{\Un\scr{V}_x}$ and $|\scr{V}_x|\leq\kappa$. Let $\scr{W}_x=\{V\in\scr{V}_x:V\meet St(x,\scr{U})\neq\es\}$ and note that $|\scr{W}_x|\leq |\scr{V}_x|\leq\kappa$.

For $x\in X$ define $W_x=\Un\scr{W}_x$. Observe that $W_x\sse St^2(x,\scr{U})$ and so $\cl{W_x}\sse\cl{St^2(x,\scr{U})}$. We show $x\in\cl{W_x}$ for each $x\in X$. For $x\in X$ let $T$ be an arbitrary open set containing $x$. There exists $U\in\scr{U}$ such that $x\in U$ and there exists $V\in\scr{V}_x$ such that $U\meet T\meet V\neq\es$. Then $V\meet St(x,\scr{U})\neq\es$ and so $V\in\scr{W}_x$. Then $T\meet W_x\neq\es$ and $x\in\cl{W_x}$.

Let $\scr{N}=\{\cl{W_x}:x\in X\}\sse\{\cl{\Un\scr{C}}:\scr{C}\in[\scr{W}]^{\leq\kappa}\}$. Then $|\scr{N}|\leq|\scr{W}|^\kappa\leq\lambda^\kappa$. As $x\in\cl{W_x}\sse\cl{St^2(x,\scr{U})}$, we see that $\scr{N}$ is a weak star-2 network for $\scr{U}$. Therefore $wsn_2(X)\leq|\scr{N}|\leq\lambda^\kappa=wL(X)^{dot(X)}$.
\end{proof}

We have as a corollary a theorem of Gotchev, Tkachenko, and Tkachuk.

\begin{corollary}[Gotchev, Tkachenko, Tkachuk~\cite{GTT16}]
If $X$ is Urysohn then $|X|\leq wL(X)^{dot(X)s\Delta_2(X)}$.
\end{corollary}

\begin{proof}
By Theorems~\ref{Xboundsn_2theta} and~\ref{somewsn2bound} we have 
$$|X|\leq wsn_2(X)^{s\Delta_2(X)}\leq\left(wL(X)^{dot(X)}\right)^{s\Delta_2(X)}=wL(X)^{dot(X)s\Delta_2(X)}.$$
\end{proof}

\begin{theorem}\label{sn3bound}
$sn_3(X)\leq wL(X)$ for any space $X$.
\end{theorem}

\begin{proof}
Let $\kappa=wL(X)$ and let $\scr{U}$ be an open cover of $X$. We want to show there exists $\scr{N}\sse\scr{P}(X)$ such that $|\scr{N}|\leq\kappa$ and for all $x\in X$ there exists $N\in\scr{N}$ such that $x\in N\sse St^3(x,\scr{U})$. As $wL(X)=\kappa$ there exists $\scr{V}\sse\scr{U}$ such that $|\scr{V}|\leq\kappa$ and $X=\cl{\Un\scr{V}}$. 

Fix $x\in X$. There exists $V_x\in\scr{V}$ such that $St(x,\scr{U})\meet V_x\neq\es$. Therefore $V_x\sse St^2(x,\scr{U})$. As $St(x,\scr{U})\meet V_x\neq\es$ there exists $U\in\scr{U}$ such that $x\in U$ and $U\meet V_x\neq\es$. It follows that $x\in U\sse St(V_x,\scr{U})$. As $V_x\sse St^2(x,\scr{U})$, we have $x\in St(V_x,\scr{U})\sse St^3(x,\scr{U})$. Thus $\scr{N}=\{St(V_x,\scr{U}):x\in X\}$ is a star-$3$ network for $\scr{U}$. We have $sn_3(X)\leq |\scr{N}|\leq |\{St(V,\scr{U}): V\in\scr{V}\}|\leq |\scr{V}|\leq\kappa=wL(X)$.
\end{proof}

Note that in the above proof the star-3 network $\scr{N}$ for $\scr{U}$ in fact consists of open sets. We have as a corollary a generalization of Proposition 4.7 in~\cite{BBR11}.

\begin{corollary}[Basile, Bella, Ridderbos~\cite{BBR11} in the case where $\Delta_3(X)=\omega$]\label{delta3}
If $X$ is a space for which $\Delta_3(X)$ is defined then $|X|\leq wL(X)^{\Delta_3(X)}$.
\end{corollary}

\begin{proof}
By Theorems~\ref{general} and~\ref{sn3bound} we have $|X|\leq sn_3(X)^{\Delta_3(X)}\leq wL(X)^{\Delta_3(X)}$.
\end{proof}

Notice that the case $m=3$ of Theorem~\ref{general} is a logical improvement of Corollary~\ref{delta3} in light of Theorem~\ref{sn3bound}.

For the remainder of this section we turn to results involving the Urysohn extent $Ue(X)$ (Definition~\ref{Uryextent}). While it is well known that $e(X)\leq L(X)$, we have a sharper upper bound for $Ue(X)$ given in the next proposition. See Definition~\ref{wLaL} for the definition of $aL(X)$.

\begin{proposition}\label{UeaL}
$Ue(X)\leq aL(X)$ for any space $X$.
\end{proposition}

\begin{proof}
Let $\kappa=aL(X)$ and suppose by way of contradiction there exists a $\theta$-closed, Urysohn discrete set $D$ such that $|D|=\kappa^+$. As $D$ is $\theta$-closed, for all $x\in X\minus D$ there exists an open set $U_x$ containing $x$ such that $\cl{U_x}\meet D=\es$. For all $x\in D$, there exists an open set $V_x$ containing $x$ such that $\cl{V_x}\meet D\minus\{x\}=\es$. Then $\{U_x:x\in X\minus D\}\un\{V_x:x\in D\}$ is an open cover of $X$. As $aL(X)=\kappa$, there exists $A\sse X\minus D$ and $B\sse D$ such that $|A|\leq\kappa$, $|B|\leq\kappa$, and $X=\Un_{x\in A}\cl{U_x}\un\Un_{x\in B}\cl{V_x}$. Pick $d\in D\minus B$. Then $d\notin\Un_{x\in A}\cl{U_x}$ by choice of each $U_x$, and $d\notin\Un_{x\in B}\cl{V_x}$ by choice of each $V_x$. As this is a contradiction, we conclude $Ue(X)\leq\kappa$. 
\end{proof}

The next example demonstrates that the spread between $Ue(X)$ and $e(X)$ can be quite large.

\begin{example}
In this example we give a Hausdorff space $X$ for which $Ue(X)<e(X)$. Let $X$ be the Kat\v{e}tov H-closed extension of the natural numbers; that is, $X=\kappa\omega$. It is well-known that this space has a closed discrete set of cardinality $2^{\mathfrak{c}}$ and therefore $e(X)=2^{\mathfrak{c}}$. However, as $X$ is H-closed we have $aL(X)=\omega$ and thus by Proposition~\ref{UeaL} we have $Ue(X)=\omega$.
\end{example}

The next theorem is a special case of the well-known Erd\H{o}s-Rado Theorem from infinite Ramsey Theory. It will be used in the proof of Theorem~\ref{Ue(X)}.

\begin{theorem}[Erd\H{o}s-Rado]\label{ER}
Let $X$ be a set, let $\kappa$ be an infinite cardinal, and let $f:[X]^2\to\kappa$ be a function. If $|X|>2^\kappa$ then there exists $Y\sse X$ and $\alpha<\kappa$ such that $|Y|=\kappa^+$ and $f(x,y)=\alpha$ for all $x\neq y$ in $Y$.
\end{theorem}

Our next result is, to the best of the author's knowledge, the first application of the Erd\H{o}s-Rado Theorem to a result involving the regular diagonal degree $\overline{\Delta}(X)$.

\begin{theorem}\label{Ue(X)}
If $X$ is Urysohn then $|X|\leq 2^{Ue(X)\overline{\Delta}(X)}$. 
\end{theorem}

\begin{proof}
Let $\kappa=Ue(X)\overline{\Delta}(X)$. For each $\alpha<\kappa$ let $U_\alpha$ be an open set in $X^2$ such that $\Delta_X=\Meet_{\alpha<\kappa}U_\alpha=\Meet_{\alpha<\kappa}\cl{U_\alpha}$. For every $\alpha<\kappa$ and $x\in X$ there exists an open set $V(x,\alpha)$ in $X$ such that $(x,x)\in V(x,\alpha)\times V(x,\alpha)\sse U_\alpha$. For all $x\neq y\in X$ there exists $\alpha(x,y)\leq\kappa$ such that $(x,y)\in X^2\minus\cl{U_{\alpha(x,y)}}$. Define the function $f:[X]^2\to\kappa$ by $f(x,y)=\alpha(x,y)$.

Suppose by way of contradiction that $|X|>2^\kappa$. Then by Theorem~\ref{ER}, there exists $\alpha<\kappa$ and $Y\sse X$ such that $|Y|=\kappa^+$ and $f(x,y)=\alpha(x,y)=\alpha$ for all $x\neq y\in Y$. We show that $Y$ is a $\theta$-closed Urysohn discrete set by showing $Y$ has no $\theta$-accumulation points. (Definition~\ref{thetaacc}). Suppose $z$ is a $\theta$-accumulation point of $Y$. Then $\cl{V(z,\alpha)}\meet Y\minus\{z\}\neq\es$ and in fact $\cl{V(z,\alpha)}\meet Y\minus\{z\}$ must contain at least two distinct points. To see this, suppose $\cl{V(z,\alpha)}\meet Y\minus\{z\}=\{p\}$. As $p\neq z$ and $X$ is Hausdorff, there exists an open set $W$ containing $z$ such that $p\notin\cl{W}$. Then $\es\neq\cl{V(z,\alpha)\meet W}\meet Y\minus\{z\}\sse\cl{V(z,\alpha)}\meet Y\minus\{z\}=\{p\},$ a contradiction since $p\notin\cl{W}$. Therefore there are two distinct points $x\neq y$ in $\cl{V(z,\alpha)}\meet Y\minus\{z\}$. It follows that $(x,y)\in\cl{V(z,\alpha)}\times\cl{V(z,\alpha)}\sse\cl{U_\alpha}$. This contradicts the fact that $(x,y)\notin\cl{U_{\alpha(x,y)}}=\cl{U_\alpha}$ as $x$ and $y$ are in $Y$. 

Therefore $Y$ has no $\theta$-accumulation points. It is straightforward to see that this implies $Y$ is $\theta$-closed and Urysohn discrete. As $|Y|=\kappa^+$, this contradicts that $Ue(X)\leq\kappa$. Therefore $|X|\leq 2^\kappa$.
\end{proof}

As $Ue(X)\leq aL(X)$ by Proposition~\ref{UeaL}, we see that Theorem~\ref{Ue(X)} is a variation on Gotchev's result from~\cite{Got19} that if $X$ is Urysohn then $|X|\leq aL(X)^{\overline{\Delta}(X)}$ (Corollary~\ref{GotaL} in the next section).

\section{On regular star networks.}

In this section we show that if $X$ is Urysohn then $|X|\leq\overline{sn}(X)^{\overline{\Delta}(X)}$ (Theorem~\ref{overlinesnCardBound}). Results of Buzyakova and Gotchev involving $\overline{\Delta}(X)$ follow as corollaries. It was shown by Buzyakova in~\cite{Buz06} that the cardinality of a ccc space with a regular $G_\delta$-diagonal has cardinality at most $\mathfrak{c}$. This was generalized by Gotchev in~\cite{Got19} who showed that if $X$ is Urysohn then $|X|\leq 2^{c(X)\overline{\Delta}(X)}$. By showing $\overline{sn}(X)\leq 2^{c(X)}$ for any space $X$ (Theorem~\ref{amazingoverlinesn}), we see that Buzyakova's result follows from Theorem~\ref{overlinesnCardBound}. Using three additional upper bounds for $\overline{sn}(X)$, three results of Gotchev also follow from Theorem~\ref{overlinesnCardBound}. Indeed, the proof of Theorem~\ref{overlinesnCardBound} contains the fundamental interaction between the cardinality of $X$ and $\overline{\Delta}(X)$ that is at the core of Buzyakova's theorem and the three theorems of Gotchev.

\begin{theorem}\label{overlinesnCardBound}
If $X$ is Urysohn then $|X|\leq\overline{sn}(X)^{\overline{\Delta}(X)}$.
\end{theorem}

\begin{proof}
Let $\kappa=\overline{\Delta}(X)$ and let $\lambda=\overline{sn}(X)$. By Proposition~\ref{gotchev} there exists a family $\{\scr{U}_\alpha:\alpha<\kappa\}$ of open covers of $X$ such that if $x\neq y\in X$ then there exists $\alpha<\kappa$ and $U\in\scr{U}_\alpha$ such that $x\in U$ and $y\notin\cl{St(U,\scr{U}_\alpha)}$.

As $\lambda=\overline{sn}(X)$, for all $\alpha<\kappa$ there exists a regular star network $\scr{N}_\alpha$ for $\scr{U}_\alpha$ such that $|\scr{N}_\alpha|\leq\lambda$. Let $\scr{N}=\Un_{\alpha<\kappa}\scr{N}_\alpha$. Then $|\scr{N}|\leq\lambda\cdot\kappa$. Let $\scr{C}=\{\Meet\scr{M}:\scr{M}\in[\scr{N}]^{\leq\kappa}\}$. Then $|\scr{C}|\leq|\scr{N}|^\kappa\leq (\lambda\cdot\kappa)^\kappa=\lambda^\kappa$. 

Fix $x\in X$. For all $\alpha<\kappa$ there exists $N_\alpha\in\scr{N}_\alpha\sse\scr{N}$ such that $x\in N_\alpha\sse\Meet\{\cl{St(U,\scr{U}_\alpha)}:x\in U\in\scr{U}_\alpha\}$. Then $x\in\Meet_{\alpha<\kappa}N_\alpha\sse\Meet_{\alpha<\kappa}\Meet\{\cl{St(U,\scr{U}_\alpha)}:x\in U\in\scr{U}_\alpha\}$. Now suppose $y\neq x$. Then there exists $\alpha<\kappa$ and $U\in\scr{U}_\alpha$ such that $x\in U$ and $y\notin\cl{St(U,\scr{U}_\alpha)}$. It follows that $y\notin\Meet_{\alpha<\kappa}\Meet\{\cl{St(U,\scr{U}_\alpha)}:x\in U\in\scr{U}_\alpha\}$ and therefore $y\notin\Meet_{\alpha<\kappa} N_\alpha$. This says $\{x\}=\Meet_{\alpha<\kappa}N_\alpha$ and as $\Meet_{\alpha<\kappa}N_\alpha\in\scr{C}$, we have $|X|\leq |\scr{C}|\leq\lambda^\kappa$.
\end{proof}

The following is Lemma 4.1 in~\cite{Got19}. It is a generalization of Lemma 2.1 in~\cite{Buz06}.

\begin{lemma}[Buzyakova~\cite{Buz06} for $X$ ccc, Gotchev~\cite{Got19}]\label{gotlemma}
Let $X$ be a space, let $\kappa=c(X)$, and let $U\times V$ be a nonempty open set in $X^2$. Let $\scr{U}$ be a collection of open boxes in $X^2$ such that $U\times V\sse\cl{\Un\scr{U}}$. Then there exists $\scr{V}=\{U_\alpha\times V_\alpha:\alpha<\kappa\}\sse\scr{U}$ such that $V\sse\cl{\Un_{\alpha<\kappa}V_\alpha}$ and $(U_\alpha\times V_\alpha)\meet (U\times V)\neq\es$ for all $\alpha<\kappa$.
\end{lemma}

\begin{theorem}[Kurepa~\cite{Kur62}]\label{kurepa}
For any space $X$, $c(X^2)\leq 2^{c(X)}$.
\end{theorem}

The proof of the next theorem represents a reduced version of the proof of Theorem 4.2 in~\cite{Got19} and is framed in a bit different way. 

\begin{theorem}\label{amazingoverlinesn}
If $X$ is any space then $\overline{sn}(X)\leq 2^{c(X)}$.
\end{theorem}

\begin{proof}
Let $\kappa=c(X)$ and let $\scr{U}$ be an open cover of $X$. Define $W=\Un\{U\times U:U\in\scr{U}\}$. By Theorem~\ref{kurepa} there exists a collection $\scr{V}$ of open boxes in $X^2\minus\cl{W}$ such that $|\scr{V}|\leq 2^\kappa$ and $X^2\minus\cl{W}\sse\cl{\Un\scr{V}}$. Define $\scr{N}=\{X\minus cl\Un\{T:S\times T\in\scr{W}\textup{ for some }S\}:\scr{W}\in[\scr{V}]^{\leq\kappa}\}$. Observe $|\scr{N}|\leq\left|[\scr{V}]^{\leq\kappa}\right|\leq\left(2^{\kappa}\right)^\kappa=2^\kappa$. We show $\scr{N}$ is a regular star network for $\scr{U}$.

Fix $x\in X$. Let $x\in U\in\scr{U}$ and set $Y=\{y\in X:(x,y)\in X^2\minus\cl{W}\}$. We show $X\minus\cl{St(U,\scr{U})}\sse Y$. Let $y\in X\minus\cl{St(U,\scr{U})}$. Then there exists an open set $R$ containing $y$ such that $R\meet St(U,\scr{U})=\es$. We show $(U\times R)\meet W=\es$. Suppose by way of contradiction that there exists $(a,b)\in (U\times R)\meet W$. Then there exists $Z\in\scr{U}$ such that $(a,b)\in (U\times R)\meet (Z\times Z)$. It follows that $a\in U\meet Z$ and $b\in R\meet Z$. Therefore $Z\sse St(U,\scr{U})$ and $b\in R\meet St(U,\scr{U})$. But $R\meet St(U\,\scr{U})=\es$ and so $(U\times R)\meet W=\es$. It follows that $(x,y)\in X^2\minus\cl{W}$, $y\in Y$, and that $X\minus\cl{St(U,\scr{U})}\sse Y$. As this holds for all $U\in\scr{U}$ such that $x\in U$, we have 
$$\{y\in X:(x,y)\in\cl{W}\}=X\minus Y\sse\Meet\{St(U,\scr{U}):x\in U\in\scr{U}\}.$$

As $\Delta_X\sse W$, we have $(x,x)\in W$. There exists an open set $B$ containing $x$ such that $B\times B\sse W$. For each $y\in Y$ there exists an open box $U_y\times V_y$ such that $x\in U_y\sse B$, $y\in V_y$, and $U_y\times V_y\sse X^2\minus\cl{W}$. As $\kappa=c(X)$, there exists $A\in [Y]^{\leq\kappa}$ such that $Y\sse\cl{\Un\{V_y:y\in A\}}$. Let $\scr{D}=\{U_y\times V_y:y\in A\}$ and note for all $y\in A$ we have $U_y\times V_y\sse X^2\minus\cl{W}\sse\cl{\Un\scr{V}}$. By Lemma~\ref{gotlemma} for each $y\in Y$ there exists $\scr{V}_y=\{U(y,\alpha)\times V(y,\alpha):\alpha<\kappa\}$ such that $\scr{V}_y\sse\scr{V}$, $V_y\sse\cl{\Un\{V(y,\alpha):\alpha<\kappa\}}$, and $(U(y,\alpha)\times V(y,\alpha))\meet (U_y\times V_y)\neq\es$. It is straightforward to see that $Y\sse\cl{\Un_{y\in A}\Un_{\alpha<\kappa}V(y,\alpha})$.

Now, for all $y\in A$ we have $U_y\sse B$ and $U_y\meet U(y,\alpha)\neq\es$ for all $\alpha<\kappa$. Thus $U(y,\alpha)\meet B\neq\es$ for all $\alpha<\kappa$. Suppose that $V(y,\alpha)\meet B\neq\es$ for some $\alpha<\kappa$. Then $(U(y,\alpha)\times V(y,\alpha)\meet (B\times B)\neq\es$. However, $B\times B\sse W$ and $U(y,\alpha)\times V(y,\alpha)\in\scr{V}$ and is a subset of $X^2\minus\cl{W}$, a contradiction. Therefore $V(y,\alpha)\meet B=\es$ for all $\alpha<\kappa$ and $y\in A$. Define $N=X\minus\cl{\Un_{y\in A}\Un_{\alpha<\kappa}V(y,\alpha})$ and note $N\in\scr{N}$. As $B\meet\cl{\Un_{y\in A}\Un_{\alpha<\kappa}V(y,\alpha})=\es$, we have $x\in N$. Also, by the last sentences of the previous two paragraphs, we have $N\sse X\minus Y\sse\Meet\{St(U,\scr{U}):x\in U\in\scr{U}\}$.

It follows that $\scr{N}$ is a regular star network for $\scr{U}$ and that $\overline{sn}(X)\leq |\scr{N}|\leq 2^\kappa=2^{c(X)}$.
\end{proof}
Observe that the regular star network in the proof of Theorem~\ref{amazingoverlinesn} in fact consists of open sets. In addition, it is interesting to note that while no diagonal degree is used in the proof, the space $X^2$, $\Delta_X$, and Kurepa's theorem are still used. 

By Theorem~\ref{overlinesnCardBound} and~\ref{amazingoverlinesn}, Gotchev's generalization of the result of Buzyakova follows.

\begin{corollary}[Buzyakova~\cite{Buz06} for $X$ ccc, Gotchev~\cite{Got19}]\label{buzyakova}
If $X$ is Urysohn then $|X|\leq 2^{c(X)\overline{\Delta}(X)}$.
\end{corollary}

\begin{proposition}\label{basicoverlinesn}
$\overline{sn}(X)\leq aL(X)$ for any space $X$.
\end{proposition}

\begin{proof}
Let $\kappa=aL(X)$ and let $\scr{U}$ be an open cover of $X$. As $\kappa=aL(X)$ there exists $\scr{V}\in [\scr{U}]^{\leq\kappa}$ such that $X=\Un_{V\in\scr{V}}\cl{V}$. We show $\scr{N}=\{\cl{V}:V\in\scr{V}\}$ is a regular star network for $\scr{U}$. Let $x\in X$ and let $U\in\scr{U}$ such that $x\in U$. There exists $V\in\scr{V}$ such that $x\in\cl{V}$. It follows that $U\meet V\neq\es$ and so $V\sse St(U,\scr{U})$. Therefore $\cl{V}\sse\cl{St(U,\scr{U})}$ and $x\in\cl{V}\sse\Meet\{\cl{St(U,\scr{U})}:x\in U\in\scr{U}\}$. Then $\scr{N}$ is a regular star network for $\scr{U}$ and since $|\scr{N}|\leq |\scr{V}|\leq\kappa$ we have $\overline{sn}(X)\leq aL(X)$.
\end{proof}

Observe that we now have three versions of $sn(X)$ that are bounded above by variations on the Lindel\"of degree: $sn(X)\leq L(X)$ (Proposition~\ref{basicsn}), $\overline{sn}(X)\leq aL(X)$ (Proposition~\ref{basicoverlinesn}), and $sn_3(X)\leq wL(X)$ (Theorem~\ref{sn3bound}). 

By Proposition~\ref{basicoverlinesn} and Theorem~\ref{overlinesnCardBound} we have the following result of Gotchev. It is the first of three results of Gotchev that follow from Theorem~\ref{overlinesnCardBound}.

\begin{corollary} [Gotchev~\cite{Got19}]\label{GotaL}
If $X$ is Urysohn then $|X|\leq aL(X)^{\overline{\Delta}(X)}$.
\end{corollary}

The proof of the following theorem represents a reduced version of the proof of Theorem 4.7 in~\cite{Got19}.

\begin{theorem}\label{overlinesnBound}
$\overline{sn}(X)\leq wL(X)^{\chi(X)}$ for any space $X$.
\end{theorem}

\begin{proof}
Let $\kappa=\chi(X)$, let $\lambda=wL(X)$, and let $\scr{U}$ be an open cover of $X$. As $\lambda=wL(X)$ there exists $\scr{W}\in[\scr{U}]^{\leq\lambda}$ such that $X=\cl{\Un\scr{W}}$. For all $x\in X$ let $\{V(x,\alpha):\alpha<\kappa\}$ be a neighborhood base at $x$. 

Fix $x\in X$. For all $\alpha<\kappa$ there exists $W(x,\alpha)\in\scr{W}$ such that $V(x,\alpha)\meet W(x,\alpha)\neq\es$. Let $\scr{W}_x=\{W(x,\alpha):\alpha<\kappa\}$ and note $|\scr{W}_x|\leq\kappa$. For all $\alpha<\kappa$ let $\scr{W}(x,\alpha)=\{W\in\scr{W}_x:W\meet V(x,\alpha)\neq\es\}$. Then for all $\alpha<\kappa$ we have $\scr{W}(x,\alpha)\neq\es$ and $|\scr{W}(x,\alpha)|\leq |\scr{W}_x|\leq\kappa$.

We show $x\in\cl{\Un\scr{W}(x,\alpha)}$. Let $B$ be an open set containing $x$. Then $x\in B\meet V(x,\alpha)$ and so there exists $\beta<\kappa$ such that $x\in V(x,\beta)\sse B\meet V(x,\alpha)$. As $W(x,\beta)\meet V(x,\beta)\neq\es$, we have $B\meet V(x,\alpha)\meet W(x,\beta)\neq\es$ and therefore $W(x,\beta)\in\scr{W}(x,\alpha)$. As $B\meet W(x,\beta)\neq\es$ we have $x\in\cl{\Un\scr{W}(x,\alpha)}$.

Now unfix $x\in X$ and let $\scr{N}=\left\{\Meet_{\alpha<\kappa}\cl{\Un\scr{M}_\alpha}:\{\scr{M}_\alpha:\alpha<\kappa\}\in\left[[\scr{W}]^{\leq\kappa}\right]^{\leq\kappa}\right\}$. Then $|\scr{N}|\leq(|\scr{W}|^\kappa)^\kappa\leq(\lambda^\kappa)^\kappa=\lambda^\kappa$. We show $\scr{N}$ is a regular star network for $\scr{U}$. First, consider $St(U,\scr{U})$ for some $x\in X$ and $U\in\scr{U}$ such that $x\in U$. Then there exists $\alpha<\kappa$ such that $V(x,\alpha)\sse U$. For $W\in\scr{W}(x,\alpha)$ we have $W\meet V(x,\alpha)\neq\es$. Therefore $W\meet U\neq\es$ which implies $W\sse St(U,\scr{U})$ and $\Un\scr{W}(x,\alpha)\sse St(U,\scr{U})$. It follows that $\cl{\Un\scr{W}(x,\alpha)}\sse\cl{St(U,\scr{U})}$ and that $x\in\Meet_{\alpha<\kappa}\cl{\Un\scr{W}(x,\alpha)}\sse\Meet\{St(U,\scr{U}):x\in U\in\scr{U}\}$. Now for all $\alpha<\kappa$ we have $\scr{W}(x,\alpha)\in[\scr{W}]^{\leq\kappa}$ and so $\Meet_{\alpha<\kappa}\cl{\Un\scr{W}(x,\alpha)}\in\scr{N}$. This shows $\scr{N}$ is a regular star network for $\scr{U}$ and that $\overline{sn}(X)\leq |\scr{N}|\leq\lambda^\kappa=wL(X)^{\chi(X)}$.
\end{proof}

By Theorems~\ref{overlinesnCardBound} and~\ref{overlinesnBound} we have a second result of Gotchev that follows from Theorem~\ref{overlinesnCardBound}.

\begin{corollary} [Gotchev~\cite{Got19}] If $X$ is Urysohn then $|X|\leq wL(X)^{\chi(X)\overline{\Delta}(X)}$.
\end{corollary}

\begin{proof}
By Theorems~\ref{overlinesnCardBound} and~\ref{overlinesnBound} we have $|X|\leq\overline{sn}(X)^{\overline{\Delta}(X)}\leq\left(wL(X)^{\chi(X)}\right)^{\overline{\Delta}(X)}=wL(X)^{\chi(X)\overline{\Delta}(X)}$.
\end{proof}

\begin{theorem}\label{anotheroverlinesnbound}
$\overline{sn}(X)\leq 2^{2^{wL(X)}}$ for any space $X$.
\end{theorem}

\begin{proof}
Let $\kappa=wL(X)$ and let $\scr{U}$ be an open cover of $X$. As $\kappa=wL(X)$ there exists $\scr{V}\sse\scr{U}$ such that $|\scr{V}|\leq\kappa$ and $X=\cl{\Un\scr{V}}$. For each $U\in\scr{U}$, let $\scr{V}_U=\{V\in\scr{V}:V\meet U\neq\es\}$. Note $\scr{V}_U\neq\es$ as $X=\cl{\Un\scr{V}}$. 

Fix $x\in X$ and suppose $x\in U\in\scr{U}$. For an open set $W$ containing $x$, we have $x\in W\meet U$ and therefore there exists $V\in\scr{V}$ such that $W\meet U\meet V\neq\es$. It follows that $V\in\scr{V}_U$. Since $W\meet V\neq\es$, we have $W\meet\Un\scr{V}_U\neq\es$. Therefore $x\in\cl{\Un\scr{V}_U}\sse\cl{St(U,\scr{U})}$. 

Define $\scr{N}=\{\Meet_{\scr{C}^\prime\in\scr{C}}\cl{\Un\scr{C}^\prime}:\scr{C}\in\scr{P}(\scr{P}(\scr{V}))\}$. Then $|\scr{N}|\leq |\scr{P}(\scr{P}(\scr{V}))|\leq 2^{2^{|\scr{V}|}}\leq 2^{2^\kappa}$. We have
$$x\in\Meet\{\cl{\Un\scr{V}_U}:x\in U\in\scr{U}\}\sse\Meet\{\cl{St(U,\scr{U})}:x\in U\in\scr{U}\}.$$
As $\Meet\{\cl{\Un\scr{V}_U}:x\in U\in\scr{U}\}\in\scr{N}$, it follows that $\scr{N}$ is a regular star network for $\scr{U}$. Therefore, $\overline{sn}(X)\leq |\scr{N}|\leq 2^{2^{wL(X)}}$.
\end{proof}

By Theorems~\ref{overlinesnCardBound} and~\ref{anotheroverlinesnbound}, we have a third result of Gotchev that follows from Theorem~\ref{overlinesnCardBound}.

\begin{corollary} [Gotchev~\cite{Got19}]
If $X$ is Urysohn then $|X|\leq 2^{\overline{\Delta}(X)2^{wL(X)}}$.
\end{corollary}

We conclude with a final observation and a question.

\begin{observation}
It was shown by Basile, Bella, Ridderbos in~\cite{BBR11} that if $X$ is Hausdorff then $|X|\leq 2^{d(X)s\Delta(X)}$. We remark that this follows immediately from a recent result of Carlson. It was shown in~\cite{C23} that if $X$ is Hausdorff then $|X|\leq 2^{d(X)w\psi_c(X)}$. (Recall the definition of $w\psi_c(X)$ in~\ref{wpsi}). 
As $w\psi_c(X)\leq\psi_c(X)\leq s\Delta(X)$, the theorem of Basile, Bella, and Ridderbos follows as a corollary. In fact, Example 2.14 in~\cite{C24} is an example of a Hausdorff space $X$ such that $2^{d(X)w\psi_c(X)}<2^{d(X)s\Delta(X)}$. This example $X$ is compact, separable, and Hausdorff with $w\psi_c(X)=\omega$ and $\psi_c(X)=\mathfrak{c}$. Therefore $s\Delta(X)\geq\psi_c(X)=\mathfrak{c}$ and $2^{d(X)w\psi_c(X)}=2^{\omega\cdot\omega}=\mathfrak{c}<2^{\mathfrak{c}}=2^{\omega\cdot\mathfrak{c}}\leq 2^{d(X)s\Delta(X)}$. This shows the improvement is a strict improvement.
\end{observation}

In connection with the above observation and noting that $we(X)\leq d(X)$ for any space $X$, the following question of Basile, Bella, and Ridderbos~\cite{BBR11} is still open.

\begin{question}[Basile, Bella, and Ridderbos~\cite{BBR11}]
If $X$ is Hausdorff is $|X|\leq 2^{we(X)s\Delta(X)}$?
\end{question}

\end{document}